# SCDA: School Compatibility Decomposition Algorithm for Solving the Multi-School Bus Routing and Scheduling Problem

Zhongxiang Wang[1], Ali Shafahi[2], Ali Haghani[3]

[1]Department of Civil and Environmental Engineering, University of Maryland - College Park, MD 20742, Email: zxwang25@umd.edu
[2]Department of Computer Science, University of Maryland - College Park, MD 20742, Email: ashafahi@cs.umd.edu
[3]Department of Civil and Environmental Engineering, University of Maryland - College Park, MD 20742. Phone: (301) 405-1963, Fax: (301) 405-2585. Email: haghani@umd.edu

**Abstract**
Safely serving the school transportation demand with the minimum number of buses is one of the highest financial goals of school transportation directors. To achieve that objective, a good and efficient way to solve the routing and scheduling problem is required. Due to the growth of the computing power, the spotlight has been shed on solving the combined problem of the school bus routing and scheduling problem. We show that an integrated multi-school bus routing and scheduling can be formulated with the help of trip compatibility. A novel decomposition algorithm is proposed to solve the integrated model. The merit of this integrated model and the decomposition method is that with the consideration of the trip compatibility, the interrelationship between the routing and scheduling sub-problems will not be lost in the process of decomposition. Results show the proposed decomposed problem could provide the solutions using the same number of buses as the integrated model in much shorter time (as little as 0.6%) and that the proposed method can save up to 26% number of buses from existing research.

# 1. Introduction

Transferring students from their homes to schools during morning trips and the opposite for afternoon trips on a daily basis is an important and expensive task. The major contributor to the transportation cost is the number of buses. In order to calculate the number of buses, two sub-problems, namely the routing and scheduling, need to be solved. To better explain the problem, we first define a few terminologies:

**A Trip:** An afternoon trip starts from a school, sequentially goes through a set of school stops associated with this school while satisfying the capacity constraint and maximum ride time constraints[1];

**A Bus:** Starts from the depot, sequentially serves a set of trips and goes back to the depot; the trips that are served by one bus are compatible with each other;

**Deadhead (between an ordered trip pair):** The travel from the last stop of the preceding trip to the first stop of the successive trip;

**Compatibility:** An ordered trip pair is compatible if the finish time (start time plus the travel time) of the preceding trip plus the deadhead between this ordered trip pair is less than or equal to the start time of the successive trip;

Given a set of stops, the distances between them, and the students at each stop, routing problem is to find a good set of trips to visit all bus stops. Then, these generated trips become inputs to the scheduling problem. The scheduling problem groups the compatible trips and serves them using the minimum number of buses. Considering the routing and scheduling problem as one joint problem, the objective for the whole system is to minimize the number of buses and total vehicle time, where the vehicle time includes the travel time of the trips and the deadhead between the trips. Due to the complexity of the integrated model, most of the literature considers the routing and scheduling as two separate problems. Because the trip compatibility is unknown in the routing stage[2], the objective for the routing problem is to minimize the number of trips (not buses) and/or total travel time. For the scheduling problem, since the travel time for the trips are fixed, the objective is to minimize the number of buses and the total deadhead. This simple decomposition disconnects vital connections between the routing and scheduling problems and produces worse solutions than the optimal solution obtained from the integrated model.

In this study, we develop a Mixed Integer Linear Programming (MILP) model for the integrated school bus routing and scheduling problem. The model is solved to optimality on small size problems to test its correctness. An advanced decomposition algorithm, namely the School Compatibility Decomposition Algorithm (SCDA), is proposed to solve the model for larger problems. SCDA is superior to the traditional decomposition methods because it considers the valuable scheduling information (the compatibility) when solving the routing problem. This 'look ahead' strategy blurs the boundary between the routing and scheduling problem and makes sure that the interrelationship between them is not lost during the process of decomposition.

The remaining of paper is structured as follows. We first present a literature review for the school bus routing problem and scheduling problem. Then we present the integrated model which simultaneously does the routing and scheduling. Later, we present the School Compatibility Decomposition Algorithm. The model and SCDA are tested on two set of problems: randomly generated problems and problems from Shafahi et al. (2017). Finally, the conclusions are presented, and some future research steps are suggested.

---

[1] The maximum ride-time constraints limit the maximum duration of each trip that includes the travel time between the stops, and the pickup and drop off time at the school and the stops.
[2] During the routing stage, the trips are being generated and therefore, their compatibilities are unknown.

## 2. Literature Review

A thorough examination of the classification of the school bus routing problem is presented in Park and Kim (2010) based on the problem characteristics and solution methods. Some of the recent work on the school bus routing and scheduling (SBRS) are listed in Table 1.

Table 1 Literature Review of school bus routing and scheduling problem

| Author | Year | Objective | Constraints | # of School | Fleet | Data |
|---|---|---|---|---|---|---|
| Fügenschuh | 2009 | NOB, TDD | SCH, BTW | M | HO | Five counties in Germany, up to 102 schools, 490 trips; Artificial, up to 10 schools and 25 trips |
| Fügenschuh | 2011 | NOB, TDD | SCH, BTW | M | HO | Same as Fügenschuh, 2009 |
| Díaz-Parra et al. | 2012 | TTD, NOT | LOG, C, MRT, SBL | S | HO | Artificial, 50 problems, each problem has 200 bus stops |
| Kim et al. | 2012 | NOB | LOG, SCH | M | HO/HT | Artificial, up to 100 schools and 562 trips for both HO and HT cases. |
| Park et al. | 2012 | NOB | LOG, C, MRT, TW | M | HT | Artificial, up to 100 schools, 2000 stops, 32048 students |
| Schittekat et al. | 2013 | TTD | LOG, C, SBL | S | HO | Artificial, up to 80 stops, 400 students |
| Caceres et al. | 2014 | TTD, NOT | LOG, C, MRT, SBL | M | HO | Williamsville Central School District, 13 schools, up to 177 stops and 1237 student per school |
| Faraj et al. | 2014 | TTD | LOG, C, MRT | M | HT | Artificial, up to 67 stops, 221 students |
| Kinable et al. | 2014 | TTD | LOG, C, MT | S | HO | Artificial, up to 40 stops and 800 students |
| Bögl et al. | 2015 | TTD, PLT | LOG, C, SBL, MWD, TSF | M | HO | Artificial, up to 8 schools, 500 students |
| Chen et al. | 2015 | NOB, TDD | LOG, SCH | S | HO/HT | Benchmark problems from Kim et al., 2012; Park et al., 2012 |
| Kang et al. | 2015 | NOB, TTD | LOG, C, MRT | M | HT | 26 students, six schools, three buses |
| Kumar and Jain | 2015 | TTD | LOG, C | S | HO | Artificial, up to 40 schools, 235 trips, 11600 students |
| Mushi et al. | 2015 | TTD | LOG, C | S | HO | 58 stops, 456 students, Dar es Salaam, Africa |
| Santana et al. | 2015 | TTD | LOG, C, TW | S | HO | 600 students, 440 nodes, Bogota, Colombia |
| Silva et al. | 2015 | TTD | LOG, C, SBL | M | HT | A city in Brazilian with 23 schools and 716 students |
| Yan et al. | 2015 | BTT, PTT, PLT | LOG, C, MNB | S | HO | Six universities (treated as one school), 400 students, Taiwan |
| Yao et al. | 2016 | TTD | LOG, C, MRT | M | HO | Artificial, up to 2 schools, 116 stops, and 1088 students |
| de Souza Lima et al. | 2017 | TTD, NOT, BAL | LOG, C | M | HT | Artificial, up to 20 schools, 150 stops; Benchmark from Park, Tae and Kim 2012 |
| Shafahi et al. | 2017 | TTD, NOT, TC | LOG, C, SBL, TCC, TW | M | HO | Artificial, up to 25 schools, 200 stops, 3656 students |
| Shafahi et al. | 2018 | NOB, PLT | SCH, MTR | M | HO | Howard county with 994 trips; A district in California with 54 trips |

**Note: Objective:** NOT: Number of trips; NOB: Number of buses; TTD: Total travel distance (or travel time); TDD: Total deadhead; BTT: Bus travel times; PTT: Passenger travel time; PLT: Penalty; BAL: Balance between each trip; TC: Trip compatibility. **Constraints**: LOG: Logistic constraints; C: Capacity Constraint; MRT: Maximum ride time; TW: Time window constraint; MWD: Maximum walk distance; SBL: Sub-tour elimination constraint; TSF: Transfer; SCH: Scheduling; BTW: school bell time window; TCC: Trip compatibility constraint; MNB: Maximum number of buses. **# of schools**: S: Single-school; M: Multi-school. **Fleet**: HO: Homogeneous fleet; HT: Heterogeneous fleet

The most widely-used objectives of the routing problem are 1) minimizing the number of trips (Bodin and Berman, 1979); 2) minimizing the total travel time (Schittekat et al., 2013; Faraj et al., 2014; Kinable et al., 2014; Mushi et al., 2015; Santana et al., 2015; Silva et al., 2015; and Yao et al., 2016); and 3) the combination of these two (Díaz-Parra et al., 2012; Caceres et al., 2014). The most prominent objective of the scheduling problem is the minimization of the number of buses and total deadhead (Fügenschuh 2009; Fügenschuh 2011; Kim et al., 2012; Chen et al., 2015). Shafahi et al. (2017) provided a new alternative to solve school bus routing problem by minimizing the number of trips and the total travel time while maximizing the trip compatibility.

The problem coverage of the recent works is summarized in Figure 1. Santana et al. (2015) and Shafahi et al. (2017) tried to solve the routing problem with time window constraints. Kim et al. (2012) showed that the scheduling problem becomes an assignment problem given the routing plan with the start times of all trips. The major effort in the literature has been made to solve the routing and scheduling problem as they are prominent problems of school bus planning. Thus, in this paper, we model the routing and scheduling problem with the assumption that the locations of the stops, the number of students at each stop, and the school bell times are known. These assumptions are valid assumptions for many real-world school bus planning problems in practice.

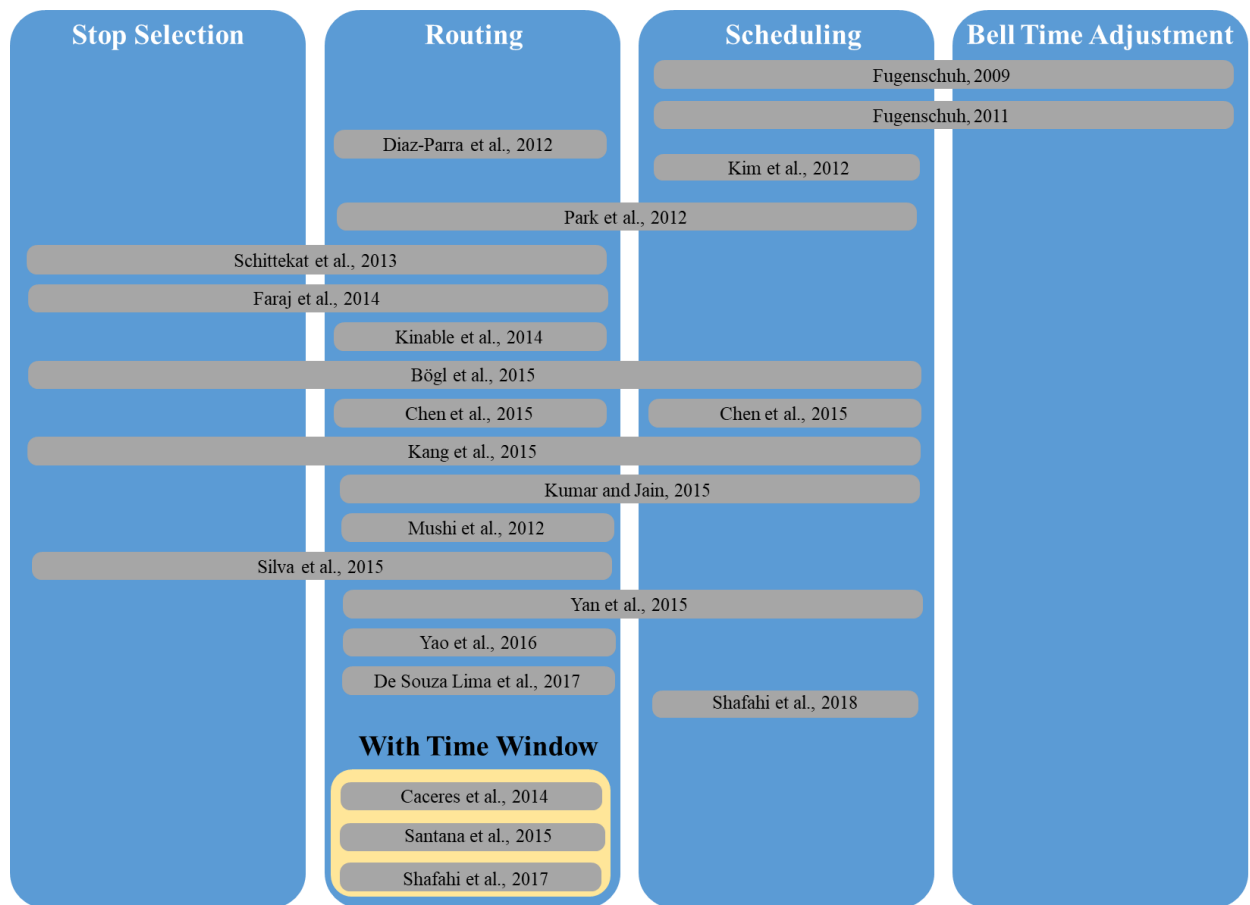

Figure 1 Problem coverage of recent school bus papers

## 3. Model Development
### 3.1 Notations

We develop a novel MILP for the integrated multi-school bus routing and scheduling problem which considers trip compatibility. The notations are summarized in *Table 2*.

Table 2 Notation summary for model formulation

**Variables for the integrated model**

| Variable | Description |
|---|---|
| $st_{s,t}$ | Binary variable, equals 1 if stop $s$ is assigned to trip $t$ |
| $ta_t$ | Binary variable, equals 1 if trip $t$ is activate (has stops assigned to it) |
| $x^t_{s1,s2}$ | Binary variable, equals 1 if in trip $t$ the bus goes directly from stop $s1$ to stop $s2$ |
| $y_{t1,t2}$ | Binary variable, equals 1 if trips $t2$ can be served after trip $t1$ (they are compatible) on the same bus |
| $y_{t,k}$ | Binary variable, equals 1 if trips $t$ is compatible with school $k$ |
| $l_{s,t}$ | Binary variable, equals 1 if the last stop of trip $t$ is stop $s$ |
| $tt_t$ | Nonnegative continuous variable, the travel time of trip $t$ |
| $dd_{t1,t2}$ | Nonnegative continuous variable, the deadhead between trip $t1$ to trip $t2$. It becomes a parameter in the scheduling problem. |
| $ddb^b_{t1,t2}$ | Nonnegative continuous variable, the deadhead between trip $t1$ to trip $t2$ on bus $b$. |
| $ddb_{t,k}$ | Nonnegative continuous variable, the travel deadhead duration from trip $t$ to school $k$. It becomes a parameter in the scheduling problem. |
| $ac^t_{s1,s2}$ | Nonnegative continuous variable, the units of "artificial commodity" that is shipped from stop $s1$ to $s2$ by trip $t$ (use for sub-tour elimination constraints) |

**Parameters for the integrated model**

| Parameter | Description |
|---|---|
| $SCH$ | Set of schools |
| $sb_k$ | The school bell (dismissal) time for school k |
| $TP[k]$ | Set of potential trips dedicated to school $k$ |
| $TP$ | Set of all potential trips |
| $SP[k]$ | Set of stops belong to school k |
| $SP$ | Set of all stops |
| $Cap$ | The capacity of each bus |
| $STU_s$ | The number of students at stop $s$ |
| $PT(stu)$ | The total pickup time for a trip with *stu* number of students |
| $DT(stu)$ | The total drop-off time for a trip with *stu* number of students |
| MNT | Minimum number of trips (for each school) |
| $AAT$ | The number of additional allowed trips (AAT) |
| $E$ | Set of compatible active trip pairs, used in Model 3 |
| $UTC[k]$ | Unassigned trip capacity for school k, used in Model 5 |
| $O\{t\}$ | School that trip t belongs to |
| $SDT$ | Start depot trip |
| $EDT$ | End depot trip |
| $D_{s1,s2}$ | The duration to drive from stop $s1$ to $s2$ |
| $M$ | A large positive value (big-M) |
| $\alpha_B$ | Coefficient for the number of buses |
| $\alpha_T$ | Coefficient for the total travel time (with students on board) |
| $\alpha_D$ | Coefficient for the total deadhead (without students on board) |
| $\alpha_N$ | Coefficient for the number of trips |
| $\alpha_C$ | Coefficient for the trip compatibilities |
| MRT | Maximum ride time |

The potential trip set ($TP[k], \forall k \in SCH$) needs to be defined first before solving the integrated problem. Although we do not know the exact cardinality of this set, we can contract the solution space by finding appropriate bounds for this unknown value. Clearly, a lower bound on the number of trips for a school can be calculated based on the bus capacity constraints. The minimum number of trips (MNT) is defined as follows:

**Minimum number of trips** (MNT) (per school): the smallest ceiling integer of the total number of students divided by the bus capacity. This quantity reflects the minimum number of trips, with respect to the bus capacities for transporting all students.

$$MNT[K] = \left\lceil \frac{\sum_{s \in SP[k]} STU_s}{Cap} \right\rceil, \forall k \in SCH. \tag{1}$$

Upper-bounding the number of trips is hard due to the maximum ride time constraint. In this paper, the upper bound is set to be MNT pluses an artificial number which we treat as a hyper-parameter – additional allowed trips (AAT).

**Additional allowed trips** (AAT): the maximum additional number of trips that can be used for each school more than the MNT (e.g. school A has 180 students and the homogeneous bus capacity is 48. Then MNT = $\lceil 180/48 \rceil$ = 4. If the AAT is set to be the same as MNT, at most eight (four plus four) trips can be used for school A). We perform sensitivity analysis on this hyper-parameter in section 5.2.

Given the lower bound (MNT) and the upper bound (MNT+AAT), the size of the potential trip set for every school is constrained ($TP[k], \forall k \in SCH$). The solution of the routing problem may not use all of the trips (MNT+AAT). After solving the routing problem, we call trips *active* if they have some stops assigned to them (the trip is a part of the solution). Other trips that have no stops assigned to them are inactive and are excluded from the scheduling problem. In addition, we ensure that the buses all start and end at the depot by adding dummy trips for the start and end depot trips and adding the appropriate constraints ($ta_{SDT} = ta_{EDT} = 1, tt_{SDT} = tt_{EDT} = 0$).

**3.2 Integrated Multi-School Bus Routing and Scheduling Model**

Our Integrated Multi-School Bus Routing and Scheduling Model (Model 1) assumes that all buses have the same capacity (homogenous fleet). Each trip only visits the stops belonging to one school (single-load) but buses can serve trips belong to different schools. Each stop should be served exactly once by one trip (single-visit). The MILP for Model 1 is as follows:

**Model 1**

$$\min z1 = \alpha_B \sum_{t \in TP} y_{SDT,t} + \alpha_T \sum_{t \in TP} tt_t + \alpha_D \sum_{t1 \in TP \cup SDT} \sum_{t2 \in TP \cup EDT \setminus t1} ddb_{t1,t2} \tag{2}$$

Subject to
(Routing constraints)

$$\sum_{j \in SP[k] \cup k \setminus s} x_{s,j}^t = \sum_{i \in SP[k] \cup k \setminus s} x_{i,s}^t, \forall k \in SCH, s \in SP[k] \cup k, \ t \in TP[k] \tag{3}$$

$$\sum_{s1 \in SP[k] \cup k} x_{s1,s}^t = st_{s,t}, \forall k \in SCH, t \in TP[k], s \in SP[k] \tag{4}$$

$$\sum_{t \in TP[k]} st_{s,t} = 1, \forall k \in SCH, s \in SP[k] \tag{5}$$

$$st_{s,t} \leq ta_t, \forall k \in SCH, s \in SP[k], t \in TP[k] \tag{6}$$

$$ta_t \leq \sum_{s \in SP[k]} st_{s,t} \ \forall k \in SCH, t \in TP[k] \tag{7}$$

$$\sum_{t \in TP[k]} \sum_{s \in SP[k]} x_{k,s}^t = \sum_{t \in TP[k]} ta_t, \forall k \in SCH \quad (8)$$

$$\sum_{s \in SP[k]} STU_s \times st_{s,t} \leq Cap, \forall k \in SCH, t \in TP[k] \quad (9)$$

$$\sum_{i \in SP[k] \cup k \setminus s} ac_{i,s}^t - \sum_{j \in SP[k] \cup k \setminus s} ac_{s,j}^t = st_{s,t}, \forall k \in SCH, t \in TP[k], s \in SP[k] \quad (10)$$

$$ac_{s1,s2}^t \leq M \times x_{s1,s2}^t, \forall k \in SCH, t \in TP[k], s1, s2 \in SP[k] \cup k \quad (11)$$

(Scheduling constraints)

$$x_{s,k}^t = l_{s,t}, \forall k \in SCH, t \in TP[k], s \in SP[k] \quad (12)$$

$$tt_t = \sum_{s1 \in SP[k] \cup k} \sum_{s2 \in SP[k] \setminus s1} x_{s1,s2}^t \times D_{s1,s2} + PT\left(\sum_{s \in SP[k]} STU_s \times st_{s,t}\right) \times ta_t + \sum_{s \in SP[k]} DT(STU_s) \times st_{s,t}, \forall k \in SCH, t \in TP[k] \quad (13)$$

$$tt_t \leq MRT, \forall t \in TP \quad (14)$$

$$dd_{t1,t2} = \sum_{s1 \in SP[O\{t1\}]} D_{s1,O\{t2\}} \times l_{s1,t1}, \forall t1 \in TP \cup SDT, t2 \in TP \cup EDT | t1 \neq t2 \quad (15)$$

$$sb_{O\{t1\}} + tt_{t1} + dd_{t1,t2} - M \times (1 - y_{t1,t2}) \leq sb_{O\{t2\}}, \forall t1, t2 \in TP | t1 \neq t2 \quad (16)$$

$$\sum_{t2 \in TP \cup EDT \setminus t1} y_{t1,t2} = ta_{t1}, \forall t1 \in TP \quad (17)$$

$$\sum_{t1 \in TP \cup SDT \setminus t2} y_{t1,t2} = ta_{t2}, \forall t2 \in TP \quad (18)$$

$$2 \times y_{t1,t2} \leq ta_{t1} + ta_{t2}, \forall t1, t2 \in TP | t1 \neq t2 \quad (19)$$

$$ddb_{t1,t2} \geq M \times (y_{t1,t2} - 1) + dd_{t1,t2}, \forall t1 \in TP \cup SDT, t2 \in TP \cup EDT | t1 \neq t2 \quad (20)$$

(Efficiency constraints)

$$tt_{t1} \times ta_{t1} \geq tt_{t2} \times ta_{t2}, \forall k \in SCH, t1, t2 \in TP[k] | t1 \leq t2 \quad (21)$$

(Domain constraints)

$$x_{s1,s2}^t \in \{0,1\}, \forall t \in TP, s1, s2 \in SP \quad (22)$$

$$st_{s,t} \in \{0,1\}, \forall s \in SP, t \in TP \quad (23)$$

$$ta_t \in \{0,1\}, \forall t \in TP \quad (24)$$

$$l_{s,t} \in \{0,1\}, \forall s \in SP, t \in TP \quad (25)$$

$$y_{t1,t2} \in \{0,1\}, \forall t1 \in TP \cup SDT, t2 \in TP \cup EDT | t1 \neq t2 \quad (26)$$

$$tt_t \geq 0, \forall t \in TP \quad (27)$$

$$dd_{t1,t2} \geq 0, \forall t1 \in TP \cup SDT, t2 \in TP \cup EDT \setminus t1 \quad (28)$$

$$ddb_{t1,t2}^b \geq 0, \forall b \in B, t1 \in TP \cup SDT, t2 \in TP \cup EDT \setminus t1 \quad (29)$$

### 3.3 Model description

The objective of Model 1 is to minimize the number of buses and total vehicle time. The number of buses equals to the number of routes that depart from the start depot trip. The total vehicle time includes the travel time of the trips (with students on board) and the deadhead (without students on board). The number of buses is of higher priority than the rest, so $\alpha_B$ is of higher magnitude than $\alpha_T$ and $\alpha_D$. In addition, the travel time with students on board is more important than the deadhead. Generally, we believe $\alpha_B \gg \alpha_T > \alpha_D$.

The constraints in Model 1 can be divided into four sets: routing constraints, scheduling constraints, efficiency constraints, and domain constraints.

The first set of constraints are the routing constraints. Constraints (3) are the conservation of flow, which guarantee that for every stop on every trip, the number of preceding stops (or the school) equals to the number of successive stops (or the school). Constraints (4) to (7) are stop-to-trip assignment constraints. Constraints (4) say that for every stop $s$ and every trip $t$, stop $s$ belongs to trip $t$ if and only if stop $s$ has a preceding stop (or school) visited on trip $t$. Constraints (5) make sure that every stop must be visited exactly once. Constraints (6) prevent the assignment of stops to the trips ($st_{s,t} = 0$) that are inactive ($ta_t = 0$). Constraints (7) deactivate a trip ($ta_t = 0$) if no stops are assigned to it ($\sum_{s \in SP[k]} st_{s,t} = 0$). Constraints (8) regulate that for each school, the number of active trips equals to the number of stops that are visited right after the school[3] (for afternoon trips). Constraints (9) are the capacity constraints, which ensure that for every trip, the sum of the students at the stops (for any given trip) is less or equal to the bus capacity. Constraints (10) and (11) are sub-tour elimination constraints that are formulated using the Artificial Commodity Flow (ACF) method (Bowerman et al., 1995).

The second part of the constraints are the scheduling constraints. First, the last stop (12) and the travel time (13) of each trip is calculated. The travel time of the trips include the travel time between stops and the pickup and drop off time for students at their home and the schools. A widely used pickup (PT at each school) and drop off (DT at each stop) time regression model was proposed by Braca et al. (1997) (travel time is calculated in seconds). It was adopted by Park et al. (2012) and Chen et al. (2015).

$$PT(stu) = 29.0 + 1.9 \times stu \tag{30}$$

$$DT(stu) = 19.0 + 2.6 \times stu \tag{31}$$

Due to the linearity of the regression, it can be directly incorporated into the model as follows:

$$tt_t = \sum_{s1 \in SP[k] \cup k} \sum_{s2 \in SP[k] \setminus s1} x^t_{s1,s2} \times D_{s1,s2} + 29 \times ta_t + (19 + 4.5 \times STU_s) \times \sum_{s \in SP[k]} st_{s,t}, \forall k \in SCH, t \in TP[k] \tag{32}$$

Constraints (13) are replaced by (32). Constraints (14) are Maximum Ride Time (MRT) constraints, which make sure that the travel time of a trip does not exceed a given duration. Constraints (15) calculate the deadhead between every pair of trips. Constraints (16) identify the compatible trip pairs ($y_{t1,t2} = 1$) if the start time (school bell time, $sb_{O\{t1\}}$) of the preceding trip plus the travel time of the preceding trip ($tt_{t1}$) and the deadhead between the two trips ($dd_{t1,t2}$) is less than or equal to the start time (school bell times, $sb_{O\{t2\}}$) of the successive trip. Constraints (17) assign exactly one successive trip to every active ($ta_{t1} = 1$) trip. Similarly, Constraints (18) assign exactly one preceding trip to every active ($ta_{t2} = 1$) trip. Constraints (19) guarantee that two trips can be served by one bus ($y_{t1,t2} = 1$) only if both trips are active. Constraints (20) is the deadhead between trips assigned to a bus. This is for the linearization,

---

[3] Without loss of generality, our formulation is for afternoon trips where the first stop after the depot is the school and students are dropped off at designated stops in the proximity of their addresses. Morning trips can be seen as the reverse of afternoon trips.

otherwise the deadhead term in the objective function would be
($\sum_{t1 \in TP \cup SDT} \sum_{t2 \in TP \cup EDT | t1 \neq t2} dd_{t1,t2} \times y_{t1,t2}$), which is a quadric term.

Constraints (21) are the efficiency constraints, which are not necessary to define the feasible region but can significantly improve the model's running time. The idea is called eliminating symmetries. Let's say a trip is generated. Labeling this trip as trip No.1 is identical to labeling this trip as trip No.2. To avoid the symmetries, constraints (21) label the longest trip (for every school) with the smallest trip ID.

Constraints (22) to (29) are variable domain constraints. We can decrease the number of variables by replacing the travel time and deadhead variables in the formulation with the right-hand sides of constraints (32) and (15), respectively. Then, constraints (32) and (15) can then be deleted. These two variables are kept to make the formulation clearer and understandable. Solvers would generally eliminate these redundant variables during their internal pre-processing phase before attempting to solve the problem using branch and bound.

## 4. School Compatibility Decomposition Algorithm

The core idea of the School Compatibility Decomposition Algorithm (SCDA) is to decompose the multi-school bus routing and scheduling problem into many compatibility-considered single-school routing problems and one multi-school scheduling problem. Many of the variables and constraints in Model 1 require information from all schools, making them indecomposable based on schools. We transform Model 1 into Model 2, which solves the same problem but is easier to decompose.

### 4.1 Model transformation

Considering that each school has a set of exclusive bus stops and trips, it is reasonable to decompose Model 1 into single-school problems. At the routing stage, the variables and constraints are fully decomposable by the schools (Equations 3 to 14). However, the scheduling problem, which aims to design the cross-school bus routes, cannot be directly decomposed based on schools. We design Compatibility Transformation (CF) to transform the cross-school variables and constraints of Model 1 into a transformed model such that the transformed model is easier to be broken down into many single-school problems. The compatibility transformed model is called the Model 2.

***Compatibility Transformation*** (CF): replace the trip-to-trip compatibility and deadhead ($y_{t1,t2}, dd_{t1,t2}, ddb_{t1,t2}$) with trip-to-school compatibility and deadhead ($y_{t,k}, dd_{t,k}, ddb_{t,k}$).

Constraints (15) to (20) need information from other schools, making them indecomposable by schools. The information required is just the departure time of the trips. Also, remember that the departure time of trips equals to the school bell time (for PM trips). If the school bell times are fixed, the trip-to-trip compatibilities can be replaced by the trip-to-school compatibilities. The trip-to-school compatibility, defined between every trip and all other schools, equals to one if the finish time of the trip (start time plus the travel time) plus the deadhead from the last stop on the trip to the school is less than or equal to the school bell time. Then constraints (15) to (20) can be replaced by the following constraints.

$dd_{t,k2} = \sum_{s \in SP[k]} D_{s,k2} \times l_{s,t}, \forall k \in SCH, t \in TP[k], k2 \in SCH \backslash k$ (33)

$sb_k + tt_t + dd_{t,k2} - M \times (1 - y_{t,k2}) \leq sb_{k2} \forall k \in SCH, t \in TP[k], k2 \in SCH \backslash k$ (34)

$\sum_{k2 \in SCH \backslash k} y_{t,k2} \leq ta_t, \forall k \in SCH, t \in TP[k]$ (35)

$$ddb_{t,k2} \geq M \times (y_{t,k2} - 1) + dd_{t,k2}, \forall k \in SCH, t \in TP[k], k2 \in SCH\backslash k \tag{36}$$

$$\sum_{t1 \in TP\backslash TP[k]} y_{t1,k} \leq \sum_{t \in TP[k]} ta_t, \forall k \in SCH \tag{37}$$

Constraints (33) calculate the deadhead between trips and other schools. Constraints (34) identify potential trip-to-school compatibilities. Constraints (35) regulate that each trip can have at most one successive school. Constraints (36) calculate the deadhead between trips to schools for a given bus. Similar to constraints (20), these constraints are needed for keeping the model linear. Constraints (37) make sure that for every school, the number of preceding trips is less than or equal to the number of active trips that this school needs. Notice that although the trip-to-school variables ($dd_{t,k}, y_{t,k}, ddb_{t,k}$) use the same notation as their trip-to-trip counterparts ($dd_{t1,t2}, y_{t1,t2}, ddb_{t1,t2}$), they are defined over a different domain and have different indices. The objective also needs to be changed accordingly. The compatibility transformed model (Model 2) is presented as follow.

**Model 2**

$$\min z2 = \alpha_N \sum_{k \in SCH} \sum_{t \in TP[k]} ta_t - \alpha_C \sum_{k \in SCH} \sum_{t \in TP[k]} \sum_{k2 \in SCH\backslash k} y_{t,k} + \alpha_T \sum_{k \in SCH} \sum_{t \in TP[k]} tt_t + \alpha_D \sum_{k \in SCH} \sum_{t \in TP[k]} \sum_{k2 \in SCH\backslash k} ddb_{t,k} \tag{38}$$

Subject to
Constraints (3) to (14) (Routing and some scheduling constraints)
Constraints (21) to (29) (Efficiency and domain constraints)
Constraints (33) to (37) (Compatibility transformed scheduling constraints)

Model 2 replaces the trip-to-trip compatibility in Model 1 with the trip-to-school compatibility. Due to the argument above, if an ordered trip pair is compatible, the preceding trip is also compatible with the school that the successive trip belongs to. The trip-to-school compatibility is identical to the trip-to-trip compatibility. Now, the question of the equality of Model 1 and Model 2 becomes whether the objective (38) is the same as objective (2). Clearly, the total vehicle time (travel time of the trips and the deadhead) is the same in two objectives. The only question is whether the number of buses can be expressed using the number of trips minus the trip-to-school compatibility.

**Lemma:** The first two components (number of trips minus the trip-to-school compatibility) in the objective (38) with equal coefficients ($\alpha_N = \alpha_C$) equals the optimal number of buses.

**Proof:** If no trips are compatible, the total number of buses equals to the number of trips ($\sum_{k \in SCH} \sum_{t \in TP[k]} y_{SDT,t} = \sum_{k \in SCH} \sum_{t \in TP[k]} ta_t$). When any two trips are compatible and blocked together, a single bus can serve both trips together. Then each realized compatible trip pair ($y_{t1,t2} = 1$) would yield in one bus saving. The total realized compatible trip pairs ($\sum_{t1 \in TP} \sum_{t2 \in TP\backslash t1} y_{t1,t2}$) is the number of buses that can be saved by blocking. Also, the trip-to-trip compatibility is identical to trip-to-school compatibility ($y_{t1,t2} = y_{t1,O\{t2\}}, \forall t1, t2 \in TP | t1 \neq t2$). Hence, the first two terms (number of trips minus the trip-to-school compatibility) in objective (38) yields the optimal number of buses.

Given the lemma and the argument above, the equality of Model 1 and Model 2 is proven. Model 2 is better than Model 1 because only one constraint (37) in Model 2 cannot be decomposed by the schools in comparison to five indecomposable constraints in Model 1.

Constraints (37) requires knowledge of the number of active trips from all schools, which is unknown in the routing stage. So, constraints (37) are dropped in the process of decomposition. The compatibility-considered single-school models can be obtained by taking the corresponding terms for one school (k ∈ SCH) in the objective and constraints of Model 2.

Due to the relaxation of constraints (37), the solutions from the decomposed single-school models may not be the optimal solution. If for every school, the number of preceding trips is less than or equal to the number of active trips this school uses, constraints (37) hold and dropping this constraint will not affect the optimality of the solution. However, it is possible to assign more preceding trips to a school with fewer active trips after dropping constraints (37). Under this circumstance, the trip-to-school compatibility is greater than the trip-to-trip compatibility. It means that the trip-to-school compatibility in the objective (38) is overestimated. Moreover, the total deadhead using trip-to-school compatibility in objective (38) is also overestimated.

## 4.2 Multi-school scheduling problem

Since the trip-to-school compatibility may be overestimated on the decomposed single-school problems, the decomposed problem is not the integrated model for the routing and scheduling problem but rather an advanced routing problem. Therefore, the scheduling problem needs to be solved after solving these decomposed problems. The trips generated by solving the decomposed single-school problems are considered as fixed input for the scheduling problem. The stop visiting sequence on each trip ($x^t_{s1,s2}$), the first and last stop on each trip ($l_{s,t}$), the travel time of each trip ($tt_t$) and deadhead between every pair of trips become known parameters. With the fixed start time of each trip, the scheduling problem (Model 3) becomes a simple assignment problem. The objective for the scheduling problem is to find the best plan to serve all the trips with the minimum number of buses and minimum total deadhead. Also, only trips belonging to the Active Trip Set (*ATP*) are considered in Model 3. ATP consists of the trips in *TP* that are active ($ta_t = 1, \forall t \in TP$). The compatible active trip pairs set $E$ is defined in the preprocessing step which is applied before solving model 3. An ordered trip pair (*t1, t2*) exists in $E$ if it satisfies: 1) *t1* and *t2* are both in ATP or the start depot trip (for *t1*) and end depot trip (for *t2*); 2) *t1* and *t2* are not the same trip nor the depot trips together; 3) trip *t1* to *t2* is compatible. The trip pairs that are not contained in $E$ will be discarded from the model.

**Model 3**

$$\min z3 = \alpha_B \sum_{t \in ATP} y_{SDT,t} + \alpha_D \sum_{(t1,t2) \in E} dd_{t1,t2} \times y_{t1,t2} \quad (39)$$

Subject to

$$\sum_{t1:(t1,t) \in E} y_{t1,t} = 1, \forall t \in ATP \quad (40)$$

$$\sum_{t2:(t,t2) \in E} y_{t,t2} = 1, \forall t \in ATP \quad (41)$$

$$y_{t1,t2} \in \{0,1\}, \forall (t1,t2) \in E \quad (42)$$

## 4.3 Model decomposition

### 4.3.1 Objective adjustment

Going back to the decomposed single-school problems. A simple approach to compensate the overestimation is to decrease the weights of the trip-to-school compatibility and total deadhead ($\alpha_C^{OA} < \alpha_C$ and $\alpha_D^{OA} < \alpha_D$) in the objective (38). This is called the ***Objective Adjustment***. The decomposed single-school problem adopting objective adjustment is Model 4. The pseudo code of the School Decomposition Algorithm (SDA) using Model 4 is Algorithm 1. This method is very similar to what Shafahi et al. (2017) did.

**Model 4**

$$\min z4 = \alpha_N \sum_{t \in TP[k]} ta_t - \alpha_C^{OA} \sum_{t \in TP[k]} \sum_{k2 \in SCH \setminus k} y_{t,k} + \alpha_T \sum_{t \in TP[k]} tt_t + \alpha_D^{OA} \sum_{t \in TP[k]} \sum_{k2 \in SCH \setminus k} ddb_{t,k} \quad (43)$$

Subject to
Constraints (3) to (14) (Routing and some scheduling constraints for school k)
Constraints (21) to (29) (Efficiency and domain constraints for school k)
Constraints (33) to (36) (Compatibility transformed scheduling constraints for school k)

**Algorithm 1**
1. Implement compatibility transformation on Model 1 to obtain Model 2
2. Drop constraints (37), decompose Model 2, and implement objective adjustment to obtain Model 4 for every school k in SCH
3. For k in SCH
   3.1 Solve Model 4 for school k
4. Report the routing plan as the combination of all active trips generated by solving Model 4 for every school k in SCH
5. Solve Model 3 given the routing plan

### 4.3.2 Compatibility assignment

We define compatibility assignment as an alternative to handle the overestimation of the trip-to-school compatibility and the total deadhead. The decomposed single-school problems adopting this idea are referred to as Model 5 ($\forall k \in SCH$).

**Compatibility Assignment** (CA): dynamically update the school Unassigned Trip Capacity (UTC) and only assign trip-to-school compatibility to the schools with positive UTC.

**Unassigned Trip Capacity** (per school): The number of trips of a school that have not been assigned to a preceding trip.

The underlying reason for overestimating trip-to-school compatibility and the total deadhead is the assignment of more preceding trips to schools with fewer active trips. Compatibility assignment will prevent this from happening with the help of UTC. When solving each single-school problem (*k*), Model 5 only assigns trips to other schools (*k2*) with positive UTC[k2]. This constraint is expressed as follows:

$$\sum_{t \in TP[k]} y_{t,k2} \leq UTC[k2], \forall k2 \in SCH \setminus k \quad (44)$$

Note that trips can only be compatible with schools with later bell times. It is desired to solve Model 5 for the school with the latest bell time first. Then, when solving Model *5* for schools with earlier bell times, we can guarantee that the number of preceding trips is less than or equal to the exact number of active trips of the successive school. In the beginning, the initial UTC is set to be the minimum number of trips (MNT, see section 3.1). Whenever one Model 5 (for school k) is solved, the UTC for all schools needs to be updated. The update rules are:

1. For the newly solved school (*k*), the updated UTC for school k is the number of active trips of that school;
   $$UTC[k] \leftarrow \sum_{t \in TP[k]} ta_t \quad (45)$$
2. For all other schools with later bell time (*k2*), the updated UTC is the old UTC minus the number of preceding trips (from the newly-solved school)
   $$UTC[k2] \leftarrow UTC[k2] - \sum_{t \in TP[k]} y_{t,k2}, \forall k2 \in SCH \setminus k \quad (46)$$

Note Model 5 is solved from the latest-dismissed school to the earliest-dismissed school. When school k's problem is newly solved, only schools with equal or later bell times are solved. Because the trips of these later dismissing schools are incompatible with school k, at this moment, no preceding trip is assigned to school k and UTC[k] is equal to the number of active trips of school k. Constraints (44) to (46) limit the assignment of trips to a school based on the succeeding school's UTC. In this case, trip-to-school compatibility and total deadhead will not be overestimated and there is no need to adjust their weights.

**Model 5**

$$\min z5 = \alpha_N \sum_{t \in TP[k]} ta_t - \alpha_C \sum_{t \in TP[k]} \sum_{k2 \in SCH \setminus k} y_{t,k} + \alpha_T \sum_{t \in TP[k]} tt_t + \alpha_D \sum_{t \in TP[k]} \sum_{k2 \in SCH \setminus k} ddb_{t,k} \quad (47)$$

Subject to
Constraints (3) to (14) (Routing and some scheduling constraints for school k)
Constraints (21) to (29) (Efficiency and domain constraints for school k)
Constraints (33) to (36) (Compatibility transformed scheduling constraints for school k)
Constraints (44) (Compatibility assignment constraints for school k)

The pseudo code of the School Compatibility Decomposition Algorithm (SDA) adopting compatibility assignment is Algorithm 2.

| Algorithm 2 |
|---|
| 1. Implement compatibility transformation on Model 1 to obtain Model 2
2. Drop constraints (37), decompose Model 2, and implement compatibility assignment to obtain Model 5 for every school k in SCH
3. Sort the "school solving sequence" based on the order of descending school bell-times
4. Initialize UTC as MNT
5. For k in sorted school solving sequence:
   5.1 Solve Model 5 for school k
   5.2 Update UTC for all schools
6. Report the routing plan as the combination of all active trips generated by solving Model 5 for every school k in SCH
7. Solve Model 3 given the routing plan |

### 4.4 Discussion of the Algorithm

A few things need to be further explained about the SCDA. First, the benefit of Model 5 over Model 4 is that Model 5 does not need to artificially adjust the weights in the objective function. The weights are usually sensitive from case to case. A bad choice of the objective adjustment could significantly diminish the solution quality. And it is usually costly to conduct the sensitivity analysis to estimate the best weight ranges – even then, there is no guarantee that the selected weights would generalize well to other problems. On the other hand, Model 5 has trip-to-school compatibility assignment constraints (44), which requires higher solution time.

Second, Model 5 prevents the assignment of more preceding trips to a school with fewer trips. Still, it does not guarantee to output an optimal solution to Model 2. One example is shown in Figure 2. Suppose at the moment we are trying to solve Model 5 for school B. The optimal solution for school C has already been obtained with 2 trips. By solving Model 5, the optimal solution for school B is two short trips where both trips are assigned as preceding trips to school C (number of trips and trip-to-school compatibility both equal to 2, with the total vehicle time of 27 minutes). In the end, when school A is solved with one trip, there is no UTC left (assume the trip from school A is not compatible with school B), three buses are needed. However, a global optimal solution to Model 2 can be obtained if we have one long trip for school B. This is not the optimal solution to Model 5 for school B since the objective for this long trip solution is worse than the two short trips solution (number of trips and trip-to-school compatibility both equal to 1 with total vehicle time of 30 minutes - see objective (47)). However, this solution will yield the best scheduling plan with only two buses.

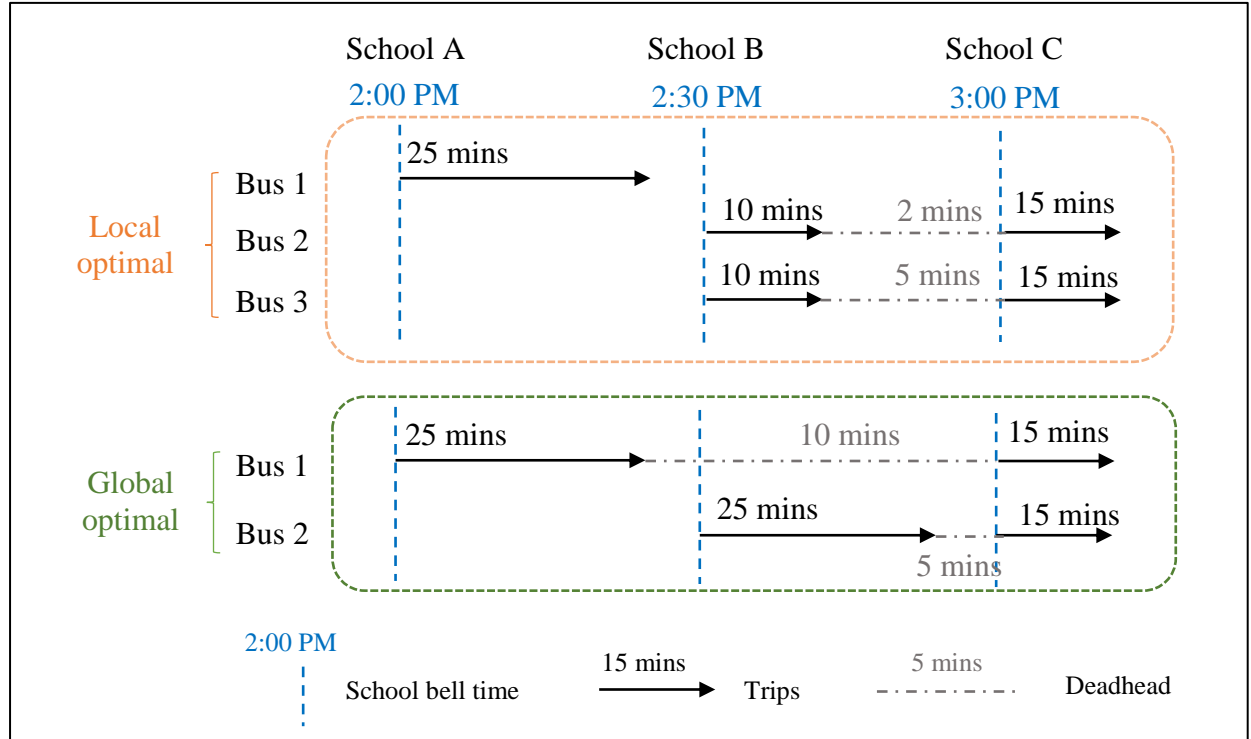

Figure 2 Illustration of Model 5

The underlying reason causing the problem just mentioned is that once a trip is assigned to a compatible successive school, there is no excess cost from adding such trip apart from its travel time. However, the trip-to-school compatibility assignment would consume a UTC, which might be used later. A simple solution is to add a small cost of adding a trip even if it has a compatible successive school. To accomplish this goal, the weights of all the components in the objective function need to be carefully tuned such that each component has different priority in the objective.

A proposal for the weights is $\alpha_B = \alpha_T = \alpha_C > \alpha_C^{CA} > \alpha_C^{OA} \gg \alpha_T > \alpha_D$ where $\alpha_C^{CA}$ is the adjusted weights for trip-to-school compatibly in Model 5. In this case, the objectives of Model 5 in the descending order of importance are: 1) minimizing the number of trips with no compatible schools (with weight $\alpha_T$); 2) minimizing the number of trips with compatible school (with weight $\alpha_T - \alpha_C^{CA}$); 3) minimizing the weighted sum of the travel time and deadhead (with weight $\alpha_T, \alpha_D$). The three objectives are of different priorities ($\alpha_T > \alpha_T - \alpha_C^{CA} \gg \alpha_T > \alpha_D$). A potential drawback for this structure is when the mentioned example does not happen: the later-dismissed school has enough active trips, the total vehicle time might increase for the purpose of reducing the number of trips. However, even under such circumstances, the number of buses will not increase.

Note Model 4 and Model 5 might both involve the weights adjustment. But the purpose is different. The weight adjustment in Model 4 is due to the overestimation of trip-to-school compatibility. But the weight adjustment in Model 5 is for the purpose of adding a small cost to the trip even if it can be assigned to a compatible school. The small cost for such trip is $\alpha_N - \alpha_C^{CA}$.

In addition, the trip-to-trip compatibility is equivalent to the trip-to-school compatibility if the bus start (departure) time equals to the school bell time. For simplicity, we do not consider the case that trips are allowed to start with a buffer after the school dismissal time. This is a realistic assumption because holding students after school involves many safety issues, which is uncommon in real-world scenarios. A slightly flexible rule would be allowing a small buffer after the school bell time. Braca et al. (1997) mentioned that the school buses should arrive at schools no later than 5 minutes after the school afternoon bell times. Fügenschuh (2009, 2011) also required that all trips arrive at school within 5 minutes of school bell times.

### 4.5 Model comparison

The models discussed above solve different versions of the problem that consider different level-of-service constraints and use a different number of variables and constraints. The comparison of these model is summarized in *Table 3*.

Table 3 Model comparison

|  |  | **Model 1** | **Model 2** | **Model 3** | **Model 4** | **Model 5** |
|---|---|---|---|---|---|---|
| | **Scope** | R, S | R, S | S | R | R |
| | **# of schools** | M | M | M | S | S |
| | **Level-of-service** | MRT | MRT | - | MRT | MRT |
| **# of variable** | **Binary** | $\|TP\|^2 + \|TP\| + \sum_{k \in SCH}(\|TP[k]\| \times \|SP[k]\|^2 + 2 \times \|TP[k]\| \times \|SP[k]\|)$ | $\|TP\| \times (\|SCH\| + 1) + \sum_{k \in SCH}(\|TP[k]\| \times \|SP[k]\|^2 + 2 \times \|TP[k]\| \times \|SP[k]\|)$ | $\|E\|$ | $\|TP[k]\| \times (\|SP[k]\|^2 + 2 \times \|SP[k]\| + \|SCH\|)$ | $\|TP[k]\| \times (\|SP[k]\|^2 + 2 \times \|SP[k]\| + \|SCH\|)$ |
| | **Continuous** | $(\|TP\| + 1)^2 + \|TP\| + \sum_{k \in SCH}(\|TP[k]\| \times \|SP[k]\|^2)$ | $2 \times \|TP\| \times (\|SCH\| + 1) + \sum_{k \in SCH}(\|TP[k]\| \times \|SP[k]\|^2)$ | 0 | $\|TP[k]\| \times (\|SP[k]\|^2 + 2 \times \|SCH\|)$ | $\|TP[k]\| \times (\|SP[k]\|^2 + 2 \times \|SCH\|)$ |
| | **# of constraints** | $2 \times \|TP\|^2 + 7 \times \|TP\| + \|SCH\| + \sum_{k \in SCH}(\|TP[k]\| \times \|SP[k]\|^2 + 0.5 \times \|TP[k]\|^2 + 5 \times \|TP[k]\| \times \|SP[k]\| + \|SP[k]\|)$ | $\|TP\| \times (3 \times \|SCH\| + 7) + \|SP\| + \sum_{k \in SCH}(\|TP[k]\| \times \|SP[k]\|^2 + \|TP[k]\|^2 + 5 \times \|TP[k]\| \times \|SP[k]\|)$ | $2 \times \|ATP\|$ | $\|TP[k]\| \times (\|TP[k]\| + \|SP[k]\|^2 + 5 \times \|SP[k]\| + 3 \times \|SCH\| + 4) + \|SP[k]\|$ | $\|TP[k]\| \times (\|TP[k]\| + \|SP[k]\|^2 + 5 \times \|SP[k]\| + 3 \times \|SCH\| + 4) + \|SP[k]\|$ |

**Note: Scope**: R (Routing), S (Scheduling); **# of school**: M (Multi-school), S (Single-school); $|set|$: the dimension of the set. Number of constraints do not include the domain (integrity, binary, nonnegative continuous) constraints
Dimension relationships:
$|SCH| \leq \sum_{k \in SCH} MNT[k] \leq |ATP| \leq \sum_{k \in SCH} MNT[k] + AAT[k] = |TP|; |ATP| \leq |SP|; |E| < 0.5 \times |ATP|^2$;
$\forall k \in SCH: MNT[k] \leq |ATP[k]| \leq MNT[k] + AAT[k] = |TP[k]|; |ATP[k]| \leq |SP[k]|;$

## 5. Computational Result
### 5.1 Experiment 1
#### *5.1.1 Traditional decomposition method*

The proposed School Compatibility Decomposition Algorithm (SCDA) differs from the traditional decomposition methods as the SCDA considers the compatibility. The traditional decomposition methods, on the other hand, simply decomposes Model 1 into routing (without compatibility, defined in Model 6) and scheduling problem (Model 3).

**Model 6**:
$$\min z6 = \sum_{t \in TP[k]} ta_t \tag{48}$$
$$\min z6' = \alpha_T \sum_{t \in TP[k]} tt_t \tag{49}$$
$$\min z6'' = \alpha_N \sum_{t \in TP[k]} ta_t + \alpha_T \sum_{t \in TP[k]} tt_t \tag{50}$$
Subject to
Constraints (2) to (10) (Routing constraints for school k)
Constraints (12) to (13) (MRT constraints for school k)
Constraints (21) to (29) (Efficiency and domain constraints for school k)

The objective (48) to (50) are possible surrogate objectives for the routing problem. The most common is minimizing the number of trips (objective 48), minimizing the total travel time (objective 49) or the combination of these two (objective 50). After solving the routing problems (with different objectives), the scheduling problems are solved using the trips generated from the routing problems. These three methods (using different routing objectives) are considered as the traditional research methods and are treated as baselines.

### 5.1.2 Experiment Setup

The usefulness of the School Compatibility Decomposition Algorithm (SCDA) depends on its solution quality and efficiency, in comparison to the integrated model and traditional decomposition methods. Thus, the first experiment is conducted to compare the solutions found by different methods in terms of the number of buses and total vehicle time on eight randomly generated dataset with increasing size. The competing methods include:

1) Mod1: Solving the Model 1 directly using commercial solver;
2) Alg1: Solving the problem using Algorithm 1 with $\alpha_C^{OA} = 5e4$;
3) Alg2: Solving the problem using Algorithm 2 without weight adjustment ($\alpha_C = 1e5$);
4) Alg2W: Solving the problem using Algorithm 2 with weight adjustment ($\alpha_C^{CA} = 9e4$);
5) MinN: Solving Model 6 and Model 3 using minimizing the number of trips as the routing objective;
6) MinTT: Solving Model 6 and Model 3 using minimizing the total travel time as the routing objective;
7) MinNT: Solving Model 6 and Model 3 using minimizing the number of trips and total travel time as the routing objective;

Eight random problems are generated with increasing size, in terms of the number of schools and stops. For each test problem, all nodes (including stops and schools) are assumed to be located within a 2-dimensional square with the length of 20 miles (105,600 ft). The locations of all nodes are designated using their longitude and latitude, which are both randomly generated as a uniformly distributed random variable between 0 and 105,600. Then, k-means algorithm is applied to group these nodes together with the number of clusters (k) equal to the number of schools. The closest nodes to the clusters' centroids are selected to be the school locations, and the rest of the nodes in each cluster are treated as bus stops for the respective school. One depot is assumed to located in the middle of this 2-dimensional square. The number of students at each stop is randomly generated as a uniformly distributed random variable between 1 and 20. The buses are assumed to have the same capacity of 66. The dismissal times for the schools are random uniform variable between 12:00 PM and 16:00 PM and they are integer multiples of 15 minutes (i.e., 2:00 PM, 3:15 PM). The bus is assumed to run at a constant speed of 20 miles per hour. The distance is Euclidean. AAT is assumed to be equal to MNT for every school.

Model 1 and all sub-problems (Model 2 to Model 6) are solved by the Gurobi Python commercial solver on a computer with Intel® Core™ i7-840 CPU, 2.93 GHz with 8 GB RAM. The code is written in Python 2.7. The parameters are $\alpha_B = \alpha_N = \alpha_C = 1e5, \alpha_T = 1, \alpha_D = 0.5, MRT = 90\ minutes$.

### 5.1.3 Results and Comparisons

The characteristics of the test problem sizes and the solutions are summarized in *Table 4*. In theory, the solution from Mod1 (if solved to optimality) should be the best among all solutions obtained by other methods. We state Mod1 found the solutions using the minimum number of buses for all scenarios and the minimum vehicle time for the scenarios if solved to optimality.

Based on the number of times that each approach found the best solution (with respect to the minimum number of buses), the approach with the highest solution quality is: Mod1 (8 times) = Alg2W (8 times) > Alg2 (7 times) > Alg1 (4 times) = MinNT (4 times) > MinTT (1 time) = MinN (1 time). In these experiments, Alg2W is the best in all decomposition approaches. It is as good as Mod1 for finding the solution with the minimum number of buses. Alg2 found solutions which required the same number of buses as what Mod1 and Alg2W require in seven scenarios out of eight. Alg1 finds the best solution (minimum number of buses) for small size problems (scenarios 1 to 4). For large problems (scenario 5 to 8), Alg1 finds solutions using at most two more buses than the best-known solution, which is better than the solutions from traditional methods (MinN, MinTT, and MinNT). The proposed method (Alg2W) can save up to 27% of buses from the best traditional method (in scenario 6, Alg2W used 8 buses as compared to 11 buses from MinNT).

In terms of the total vehicle time (sum across all scenarios), the approaches are ranked in the order of: Mod1 (6924 min) < Alg2W (6951 min) < Alg1 (6986 min) < Alg2 (7010 min) < MinNT (7357 min) < MinTT (7623 min) < MinN (10845 min). The total vehicle time for the integrated model (Mod1) is the smallest among all approaches if the optimality gap is small (scenario 1 to 7). But in scenario 8 where the optimality gap for Mod1 is relatively higher, the proposed methods (Alg1, Alg2, and Alg2W) are better than Mod1 in terms of the total vehicle time. Except for the Mod1, the list of approaches with the lowest total vehicle time is very similar to the list with the minimum number of buses. It means the proposed models (Alg1, Alg2, and Alg2W) outperform the traditional methods in terms of all criteria (the number of buses and total vehicle time). The objective of MinNT shares the same terms of Alg1 as minimizing the number of trips and total vehicle time except for the trip-to-school compatibility. And the result in *Table 4* shows that Alg1 always finds better or equal solutions (with a fewer number of buses and/or shorter vehicle time) than what MinNT does, which demonstrates the importance of compatibility in the routing problems.

Table 4 Result summary for experiment 1

| Scenario | | 1 | 2 | 3 | 4 | 5 | 6 | 7 | 8 |
|---|---|---|---|---|---|---|---|---|---|
| # of schools | | 2 | 4 | 6 | 8 | 10 | 15 | 20 | 30 |
| # of stops | | 20 | 40 | 60 | 80 | 100 | 150 | 200 | 300 |
| Mod 1 | NOB | 4* | 4* | 7* | 9* | 8* | 8* | 11* | 16* |
| | NOT | 4 | 8 | 11 | 18 | 26 | 35 | 44 | 67 |
| | TVT | **310** | **437** | **583** | **717** | **814** | **1089** | **1259** | **1715** |
| | RT | 1.27 | 18.73 | 1192 | 2585 | 1585 | 23795 | 5697 | 36000 |
| | Gap | 0 | 0 | 0 | 0.05 | 0.03 | 0.02 | 0.01 | 12.31 |
| Min N | NOB | 4* | 5 | 8 | 11 | 10 | 12 | 17 | 22 |
| | NOT | 4* | 8 | 11 | 18 | 21 | 30 | 41 | 60 |
| | TVT | **354** | **703** | **971** | **1212** | **1238** | **1638** | **2094** | **2635** |
| | RT | 1.21 | 1.87 | 2.92 | 7.25 | 10.09 | 1254 | 13.74 | 13.98 |
| Min TT | NOB | 5 | 6 | 12 | 9* | 11 | 13 | 15 | 21 |
| | NOT | 5 | 12 | 18 | 21 | 26 | 44 | 55 | 82 |
| | TVT | **330** | **503** | **707** | **733** | **877** | **1174** | **1381** | **1918** |
| | RT | 1.11 | 1.91 | 5.91 | 12.15 | 27.52 | 23.59 | 43.19 | 33.42 |
| Min NT | NOB | 4* | 4* | 7* | 9* | 10 | 11 | 14 | 18 |
| | NOT | 4* | 8 | 11 | 18 | 21 | 30 | 41 | 60 |
| | TVT | **313** | **474** | **652** | **773** | **859** | **1136** | **1391** | **1759** |
| | RT | 1.34 | 4.56 | 5.98 | 9.8 | 14.85 | 18.22 | 31.68 | 28.38 |
| Alg 1 | NOB | 4* | 4* | 7* | 9* | 10 | 9 | 13 | 18 |
| | NOT | 4 | 8 | 11 | 18 | 21 | 30 | 41 | 60 |
| | TVT | **332** | **467** | **604** | **730** | **819** | **1056** | **1327** | **1651** |
| | RT | 1.07 | 3.39 | 12.65 | 38.11 | 41.71 | 50.57 | 60.18 | 83.43 |
| Alg 2 | NOB | 4* | 4* | 7* | 9* | 8* | 8* | 13 | 16* |
| | NOT | 4 | 8 | 11 | 18 | 22 | 35 | 43 | 62 |
| | TVT | **332** | **451** | **604** | **732** | **813** | **1092** | **1296** | **1690** |
| | RT | 1.11 | 4.61 | 10.06 | 48.92 | 193.37 | 58.91 | 60.15 | 85.01 |
| Alg 2W | NOB | 4* | 4* | 7* | 9* | 8* | 8* | 11* | 16* |
| | NOT | 4 | 8 | 11 | 18 | 22 | 34 | 44 | 61 |
| | TVT | **313** | **465** | **589** | **732** | **816** | **1087** | **1290** | **1659** |
| | RT | 1.23 | 5.15 | 9.77 | 137.97 | 543.58 | 71.96 | 177.92 | 208.57 |

**Note:** * minimum number of buses; Gap: Optimal gap (%); NOB: Number of buses; NOT: Number of trips
TVT: Total vehicle time (minutes); RT: Running time (second)

  The comparison between Alg1, Alg2, and Alg2W is clear. Alg2W finds the best solutions among the three approaches with respect to the number of buses and total vehicle time for all scenarios. Alg2 is the second and Alg1 is the worst of the three. But the difference between the solutions is not significant, especially for small problems. Thanks to the additional compatibility assignment, Alg2 outperforms the Alg1. Furthermore, the weight adjustment in Alg2W improves the model and helps in finding better solutions than Alg2. Alg2W takes more time to solve than Alg2 and Alg1. The small running time increase is negligible considering the huge financial benefits of saving even just one bus. Overall, Alg2W is the best. Also, solving Alg2W only takes less than 0.57% of the time required to solve Model 1 directly (at scenario 8 where Model 1 is not even solved to optimality).

### *5.1.4 Sensitivity analysis*

We mentioned that the difference between Alg1, Alg2, and Alg2W can be from the choice of the adjusted weight in the objective of Alg1 (see section 4.3.1). The choice of the trip-to-school compatibility and deadhead weights could impact the solution quality. To illustrate this, we conduct sensitivity analysis on these weights for scenarios 7 and 8. The weight of the number of trips and trip-to-school compatibility is of higher magnitude than the travel time and deadhead. Therefore, this sensitivity analysis concentrates on the adjusted weight of the trip-to-school compatibility. The weight of the number of trips, the total travel time and deadhead are fixed ($\alpha_N = 1e5, \alpha_T = 1, \alpha_D = 0.5$). The result is presented in *Table 5*. The result shows that solutions found with different weights are very similar with respect to the number of buses, number of trips and the total vehicle time. It means that Model 4 and Algorithm 1 are not extremely sensitive to the adjusted weight. Still, different solutions are obtained under different weights for trip-to-school compatibility. And there is no clear generic trend of the adjusted weights range among all scearnios. Thus, the adjusted weight should be determined based on the specific problem.

Table 5 Sensitivity analysis for objective adjustment for Model 4

| $\alpha_C^{OA}$ | Scenario 7 | | | | Scenario 8 | | | |
|---|---|---|---|---|---|---|---|---|
| | **NOB** | **NT** | **TVT** | **RT** | **NOB** | **NT** | **TVT** | **RT** |
| 1 | **15** | 41 | 391 | 46.47 | **18** | 60 | 1776 | 39.41 |
| 10 | **15** | 41 | 1369 | 51.11 | **18** | 60 | 1698 | 40.13 |
| 100 | **14** | 41 | 1337 | 52.42 | **18** | 60 | 1722 | 43.21 |
| 1000 | **13** | 41 | 1376 | 40.33 | **18** | 60 | 1710 | 46.11 |
| 10000 | **13** | 41 | 1319 | 36.15 | **18** | 60 | 1745 | 48.38 |
| 20000 | **13** | 41 | 1322 | 36.80 | **18** | 60 | 1726 | 55.17 |
| 30000 | **13** | 41 | 1295 | 61.92 | **18** | 60 | 1693 | 79.13 |
| 40000 | **13** | 41 | 1328 | 70.72 | **18** | 60 | 1667 | 89.00 |
| 50000 | **13** | 41 | 1327 | 60.18 | **18** | 60 | 1651 | 83.43 |
| 60000 | **13** | 42 | 1332 | 82.09 | **18** | 61 | 1698 | 146.81 |
| 70000 | **13** | 42 | 1303 | 105.62 | **18** | 61 | 1676 | 97.52 |
| 80000 | **13** | 42 | 1303 | 532.15 | **17** | 61 | 1692 | 213.54 |
| 90000 | **13** | 42 | 1314 | 143.91 | **17** | 61 | 1688 | 295.00 |

**Note:** NOB: Number of buses; NOT: Number of trips; TVT: Total vehicle time (minutes); RT: Running time (second)

## 5.2 Experiment 2
### *5.2.1 Experiment Setup*

Now, the question is whether the conclusion from experiment 1 is valid for other problems and whether different parameters will affect the performance of Algorithm 2 with weight adjustment. In this set of experiments, the same eight set of mid-size problems generated in Shafahi et al. (2017)[4] (abbreviated hereon as Shafahi) is used to test the performance of the School

---
[4] To make the comparison as close to Shafahi et al. (2017) as possible, some of the tests ignore the MRT constraints because it was not considered in Shafahi et al. (2017).

Compatibility Decomposition Algorithm (SCDA). Since Algorithm 2 with weight adjustment is the best among the three proposed methods, only Algorithm 2 with weight adjustment is used in this experiment. AAT is also set to vary between zero and two. Other parameters are the same as in experiment 1. Eight combinations of different parameters settings are used:

1) A0TL15: zero AAT per school and 15 seconds running time limit per Model 5;
2) A0TL30: zero AAT per school and 30 seconds running time limit per Model 5;
3) A1TL15: one AAT per school and 15 seconds running time limit per Model 5;
4) A1TL30: one AAT per school and 30 seconds running time limit per Model 5;
5) A1TL120: one AAT per school and 120 seconds running time limit per Model 5;
6) A1TL30MRT: one AAT per school, 30 seconds running time limit per Model 5 with MRT=40 minutes;
7) A2TL30MRT: two AAT per school, 30 seconds running time limit per Model 5 with MRT=40 minutes;
8) A3TL30MRT: three AAT per school, 30 seconds running time limit per Model 5 with MRT=40 minutes;

In the approach names, A is the additional allowed trips (AAT) per school, TL is the maximum running time limit (second) for each single-school problem, and MRT is the maximum ride time parameter (for all trips) set to 40 minutes. All the case studies were run on a computer with i7 CPU 870 @ 2.93 GHz and 8GB RAM. The commercial solver used to solve single school routing and scheduling problem is FICO Xpress. The configurations and running times for these test problems are shown in Table 6.

Table 6 Configurations and running time comparison between Shafahi and SCDA

| Scenario | # of schools | # of stops | ① | ② | ③ | Running time (min) | | |
|---|---|---|---|---|---|---|---|---|
| | | | | | | Shafahi | SCDA | Ratio[④] |
| 1 | 20 | 100 | 91.4 | 13 | 0-30 | 0.28 | 4.12 | 14.71 |
| 2 | 20 | 200 | 89.6 | 16 | 0-30 | 3.46 | 9.31 | 28.70 |
| 3 | 20 | 100 | 120.7 | 13 | 0-30 | 5.52 | 23.85 | 4.32 |
| 4 | 20 | 100 | 182.8 | 13 | 0-30 | 25 | 72.90 | 2.92 |
| 5 | 25 | 125 | 90.4 | 13 | 0-30 | 3.5 | 11.93 | 3.41 |
| 6 | 20 | 100 | 91.6 | 13 | 0-90 | 0.29 | 9.52 | 32.81 |
| 7 | 20 | 200 | 89.5 | 16 | 0-90 | 17.4 | 65.55 | 3.77 |
| 8 | 20 | 200 | 91.1 | 14 | 0-16 | 4.18 | 96.10 | 22.99 |

Note: ① Average number of student per school (students per school); ② Maximum number of stops to each school (stops per school); ③ School dismissal time range (min); ④ Ratio = running time (SCDA) / running time (Shafahi)

SCDA is slower than Shafahi because SCDA dynamic updates the UTC and prevents the assignment of more preceding trips to a school with fewer trips. However, the run-time is still acceptable (less than 100 minutes) given the planning nature of the problem.

*5.2.2 Solution Quality*
The results of the solutions found by different methods are listed in Table 7, including the minimum number of buses, the mean travel time per trip and the maximum travel time per trip.

Table 7 Computational result for experiment 2

| Scenario | | 1 | 2 | 3 | 4 | 5 | 6 | 7 | 8 |
|---|---|---|---|---|---|---|---|---|---|
| MinN | NOB | **38** | **38** | **45** | **60** | **41** | **21** | **26** | **45*** |
| | AvgTT | 30.42 | 24.28 | 22.53 | 21.08 | 37.88 | 30.42 | 24.28 | 136.23 |
| | MaxTT | 79.73 | 61.89 | 58.39 | 69.34 | 100.58 | 79.73 | 61.89 | 243.69 |
| MinTT | NOB | **32** | **32** | **38** | **52** | **35** | **17** | **24** | **45*** |
| | AvgTT | 14.11 | 7.08 | 11.72 | 9.10 | 18.57 | 14.13 | 11.17 | 35.15 |
| | MaxTT | 35.33 | 15.95 | 33.07 | 29.56 | 39.62 | 35.33 | 17.26 | 64.65 |
| Shafahi | NOB | **31** | **28** | **39** | **46** | **33** | **17*** | **18** | **45*** |
| | AvgTT | 14.48 | 11.38 | 11.85 | 9.58 | 19.49 | 14.61 | 11.27 | 35.43 |
| | MaxTT | 33.07 | 19.99 | 33.07 | 29.56 | 60.51 | 33.07 | 19.08 | 68.42 |
| A0TL15 | NOB | **23*** | **26*** | **36** | **45** | **33** | **17*** | **16*** | **45*** |
| | AvgTT | 18.84 | 11.30 | 11.38 | 9.20 | 19.12 | 14.27 | 11.07 | 35.45 |
| | MaxTT | 53.11 | 24.56 | 33.07 | 29.56 | 60.51 | 33.07 | 21.53 | 68.42 |
| A0TL30 | NOB | **23*** | **26*** | **36** | **44** | **33** | **17*** | **16*** | **45*** |
| | AvgTT | 18.84 | 11.49 | 11.41 | 9.23 | 19.12 | 14.27 | 11.10 | 35.45 |
| | MaxTT | 53.11 | 25.50 | 33.07 | 29.56 | 60.51 | 33.07 | 22.77 | 68.42 |
| A1TL15 | NOB | **23*** | **26*** | **34** | **43** | **31*** | **17*** | **17** | **45*** |
| | AvgTT | 19.16 | 10.74 | 10.73 | 9.38 | 18.58 | 14.27 | 10.94 | 35.48 |
| | MaxTT | 53.11 | 24.72 | 33.07 | 29.56 | 60.51 | 33.07 | 24.27 | 68.42 |
| A1TL30 | NOB | **23*** | **26*** | **35** | **42*** | **31*** | **17*** | **18** | **45*** |
| | AvgTT | 18.01 | 10.60 | 11.01 | 8.96 | 17.35 | 14.33 | 10.45 | 35.44 |
| | MaxTT | 53.11 | 23.71 | 33.07 | 29.56 | 39.62 | 33.07 | 17.75 | 68.42 |
| A1TL120 | NOB | **23*** | **27** | **35** | **43** | **31*** | **17*** | **17** | **45*** |
| | AvgTT | 18.18 | 10.64 | 11.14 | 9.27 | 18.17 | 14.27 | 10.84 | 35.43 |
| | MaxTT | 53.11 | 20.10 | 33.07 | 29.56 | 63.95 | 33.07 | 20.14 | 68.42 |
| A1TL30 MRT | NOB | **28** | **26*** | **35** | **42*** | **31*** | **17*** | **17** | **50** |
| | AvgTT | 13.41 | 10.45 | 11.30 | 8.74 | 17.76 | 14.34 | 10.92 | 31.40 |
| | MaxTT | 33.07 | 22.67 | 33.07 | 29.56 | 39.62 | 33.07 | 19.07 | 39.97 |
| A2TL30 MRT | NOB | **28** | **26*** | **33*** | **43** | **31*** | **17*** | **18** | **50** |
| | AvgTT | 13.43 | 10.66 | 10.98 | 8.65 | 17.76 | 14.35 | 10.88 | 31.40 |
| | MaxTT | 33.07 | 24.75 | 33.07 | 29.56 | 39.62 | 33.07 | 19.21 | 39.97 |
| A3TL30 MRT | NOB | **27** | **26*** | **33*** | **43** | **31*** | **17*** | **17** | **50** |
| | AvgTT | 12.60 | 10.64 | 10.67 | 9.03 | 17.76 | 14.41 | 11.01 | 31.40 |
| | MaxTT | 33.07 | 23.42 | 33.07 | 29.92 | 39.62 | 33.07 | 20.85 | 39.97 |
| Note: **MinN**: Solution of bus scheduling problem given the routing trips, which is obtained by minimizing number of trips; **MinTT**: Solution of bus scheduling problem given the routing trips, which is obtained by minimizing aggregated travel time; **Shafahi**: Number of buses from approach MaxCom+TT(AS) in Shafahi et al. (2017); **AvgTT**: average travel time per trip (minutes); **MaxTT**: maximum travel time per trip (minutes); **\***: Minimum number of buses among all approaches for each scenario | | | | | | | | | |

The result shows that SCDA always finds better or at least equal solutions to the Shafahi test examples with respect to the number of buses for all test scenarios. The highest saving occurs at scenario 1 where 8 buses can be saved out of 31, which is equivalent to 26% improvement of the solution. The mean travel time for SCDA is usually small, less than 20 minutes. The underlying reason is that shorter trips are easier to be compatible with other trips and that minimizing total travel time is also in the objective. Therefore, SCDA prefers forming short trips, which is a good practice for school bus problem.

The results in Table 7 also reveal the tradeoff between the travel time increase and the bus saving. In scenario 1, A1TL30 found solutions with 23 buses compared to the 31 buses from Shafahi at the expense of the 3.5 minutes and 20 minutes increase of the mean and maximum travel time respectively. When limiting the maximum ride time (from A1TL30MRT), average travel time is 0.7 minutes shorter than that from Shafahi and the maximum travel time is the same as the solutions from Shafahi, but the number of buses increases to 28. Considering the high cost of a school bus and a driver and relatively low cost for a small increase in travel time, from a financial point of view, the savings gained by using fewer buses could easily justify the additional travel times. There are many scenarios that SCDA finds solutions that use fewer buses and also have a smaller mean and maximum ride time. It shows the SCDA can find better results than Shafahi with respect to all criteria (fewer number of buses and shorter mean and maximum ride time).

A cross-examination of the result from Shafahi and that from A0TL15 with respect to the bus saving, mean and maximum travel time (per trip) increase percentage is shown in FIGURE 5. Solutions from SCDA have shorter mean travel time in scenarios 2 to 7 than those from Shafahi. But on the contrary, at scenarios (like 1 and 7) where SCDA can save a significant number of buses, the maximum travel time increases. This again shows the tradeoff between bus saving and travel time increase.

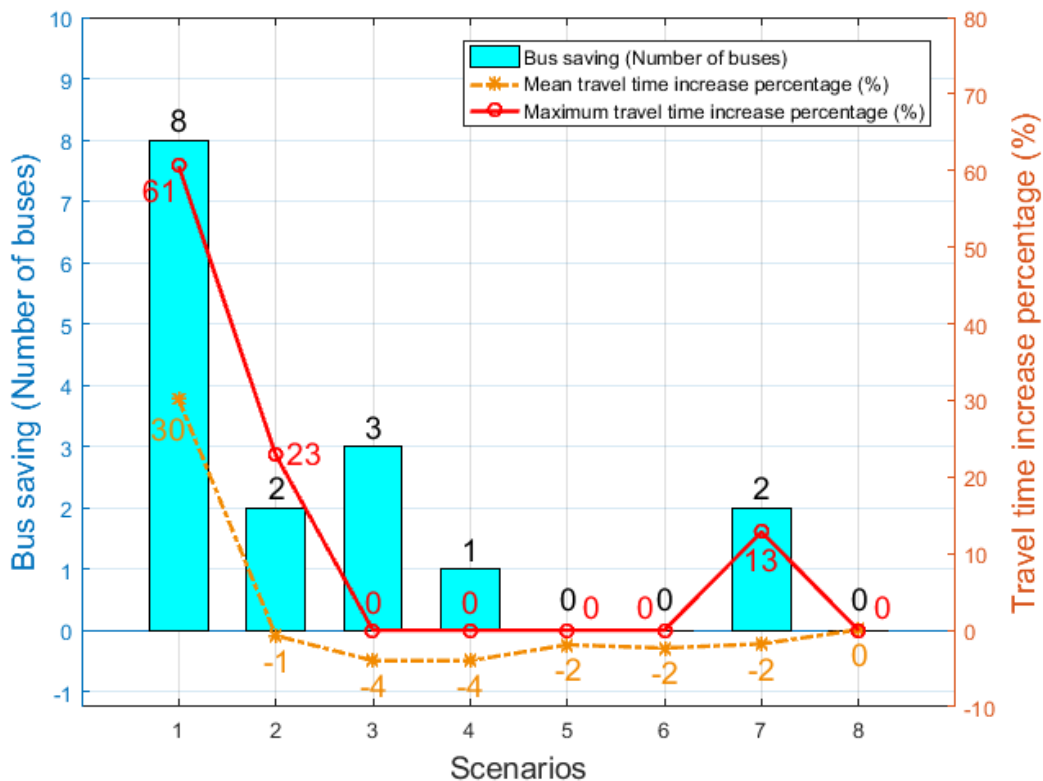

Figure 3 Comparison between solution from Shafahi and SCDA (A0TL15)

*5.2.3 Maximum Ride Time Constraint*

It can be seen that thanks to minimizing the travel time in the objective, many of the maximum travel time of the trips are still under 40 minutes even without the maximum ride time (MRT) constraint (Table 7). However, there are certain situations in which merely minimizing total travel time in the objective is not enough. Thus, we need to incorporate MRT constraint. For example, in scenario 8, without MRT, 45 buses can accommodate the school transportation demand with the highest 68 minutes per trip. By adding the MRT constraint, the maximum travel time per trip is significantly reduced (under 40 minutes), but the number of buses increases to 50.

# 6. Conclusion

In this paper, an Integrated Multi-school bus routing and scheduling model is proposed. A School Compatibility Decomposition Algorithm is developed to solve the integrated model with the consideration of trip compatibility. The biggest contribution of the model and algorithm is that the interrelation between the routing and scheduling is kept even in the decomposed problems. The validity of the model and the efficiency of the SCDA algorithm are tested on the randomly generated problems and a set of test problems developed by Shafahi et al. (2017). The first experiments show that SCDA can find solutions as good as the integrated model (in terms of the number of buses) in much shorter time (as little as 0.6%) and that it also outperforms the traditional decomposition algorithms. The second experiments show that the SCDA can find better results than Shafahi et al. (2017) with a fewer number of buses (up to 26%), and shorter mean and maximum travel time per trip (up to 7%). A few directions for future work can be identified. One of them is a more efficient algorithm to solve each single school routing problem such that it can handle more complicated problems with more stops to every single school. Another one is that a more flexible way to handle bus service start time can be devised, especially for morning trips. An appropriate time window might be more financially beneficial than a fixed service start time.

**REFERENCES**
1. Bodin, L.D., Berman, L., 1979. Routing and scheduling of school buses by computer. Transportation Science 13 (2), 113–129.
2. Bögl, M., Doerner, K.F., and Parragh, S.N., 2015. The school bus routing and scheduling problem with transfers. Networks, 65(2), pp.180-203.
3. Bowerman, R., Hall, B., and Calamai, P., 1995. A multi-objective optimization approach to urban school bus routing: Formulation and solution method. Transportation Research Part A: Policy and Practice, 29(2), pp.107-123.
4. Braca, J., Bramel, J., Posner, B. and Simchi-Levi, D., 1997. A computerized approach to the New York Cityschool bus routing problem. IIE transactions, 29(8), pp.693-702.
5. Caceres, H., Batta, R., and He, Q. "School Bus Routing with Stochastic Demand and Duration Constraints." submitted to Transportation Science (2014).
6. Chen, X., et al., 2015. Exact and metaheuristic approaches for a bi-objective school bus scheduling problem. PloS one, 10(7), p.e0132600.
7. de Souza Lima, F.M., et al., 2017. A multi-objective capacitated rural school bus routing problem with heterogeneous fleet and mixed loads. 4OR, pp.1-28.
8. Díaz-Parra, O., et al., 2012, November. A vertical transfer algorithm for the school bus routing problem. In Nature and Biologically Inspired Computing (NaBIC), 2012 Fourth World Congress on (pp. 66-71). IEEE.

9. Faraj, M.F., et al., 2014, October. A real geographical application for the school bus routing problem. In Intelligent Transportation Systems (ITSC), 2014 IEEE 17th International Conference on (pp. 2762-2767). IEEE.
10. Fügenschuh, A., 2009. Solving a school bus scheduling problem with integer programming. European Journal of Operational Research, 193(3), pp.867-884.
11. Fügenschuh, A., 2011. A set partitioning reformulation of a school bus scheduling problem. Journal of Scheduling, 14(4), pp.307-318.
12. Kang, M., et al., 2015. Development of a genetic algorithm for the school bus routing problem. International Journal of Software Engineering and Its Applications, 9(5), pp.107-126.
13. Kim, B.I., Kim, S. and Park, J., 2012. A school bus scheduling problem. European Journal of Operational Research, 218(2), pp.577-585.
14. Kinable, J., Spieksma, F.C., and Vanden Berghe, G., 2014. School bus routing—a column generation approach. International Transactions in Operational Research, 21(3), pp.453-478.
15. Kumar, Y. and Jain, S., 2015, September. School bus routing based on branch and bound approach. In Computer, Communication, and Control (IC4), 2015 International Conference on (pp. 1-4). IEEE.
16. Mushi, A.R., Mujuni, E. and Ngonyani, B., 2015. Optimizing Schedules for School Bus Routing Problem: the case of Dar Es Salaam Schools. International Journal of Advanced Research in Computer Science, 6(1).
17. Park, J., and Kim, B.I., 2010. The school bus routing problem: A review. European Journal of Operational Research, 202(2), pp.311-319.
18. Park, J., Tae, H. and Kim, B.-I., 2012. A post-improvement procedure for the mixed load school bus routing problem. European Journal of Operational Research, 217(1), pp.204-213.
19. Santana, L., Ramiro, E. and Romero Carvajal, J.D.J., 2015. A hybrid column generation and clustering approach to the school bus routing problem with time windows. Ingeniería, 20(1), pp.101-117.
20. Schittekat, P., Kinable, J., Sörensen, K., Sevaux, M., Spieksma, F. and Springael, J., 2013. A metaheuristic for the school bus routing problem with bus stop selection. European Journal of Operational Research, 229(2), pp.518-528.
21. Shafahi, A., Aliari, S., and Haghani, A., 2018. Balanced scheduling of school bus trips using a perfect matching heuristic. Transportation Research Board 97th Annual Meeting, Jan 2018, in Washington, D.C. arXiv preprint arXiv:1708.09338.
22. Shafahi, A., Wang, Z. and Haghani, A., 2017. Solving the school bus routing problem by maximizing trip compatibility. Transportation Research Record: Journal of the Transportation Research Board, (2667), pp.17-27. DOI: 10.3141/2667-03.
23. Silva, C.M., et al., 2015, September. A Mixed Load Solution for the Rural School Bus Routing Problem. In Intelligent Transportation Systems (ITSC), 2015 IEEE 18th International Conference on (pp. 1940-1945). IEEE.
24. Yan, S., Hsiao, F.Y. and Chen, Y.C., 2015. Inter-school bus scheduling under stochastic travel times. Networks and Spatial Economics, 15(4), pp.1049-1074.
25. Yao, B., Cao, Q., Wang, Z., Hu, P., Zhang, M. and Yu, B., 2016. A two-stage heuristic algorithm for the school bus routing problem with mixed load plan. Transportation Letters, 8(4), pp.205-219.